\documentclass{jcmlatex}

\begin{document}

\markboth{}{}

\title{New discovery of a conversion of generalized singular values}

\author{Weiwei Xu \thanks{Nanjing University of information Science and Technology, Nanjing, People Republic of China.} \and Yingzhou Lee}

\maketitle

\begin{abstract}
In this paper, we give a new discovery of conversion of generalized singular values. New formulation is proposed.
\end{abstract}

\section{Introduction}
The generalized singular value decomposition (GSVD) used in mathematics and numerical computations is a very useful and versatile tool. The GSVD of two matrices having the
same number of columns was first proposed by Van Loan. It is very useful in many matrix
computation problems and practical applications, such as the Kronecker canonical form of a
general matrix pencil, the linearly constrained least-squares problem, the general Gauss-Markov linear model, the generalized total least-squares problem, real-time signal processing, comparative analysis of DNA microarrays and so on[1,2,3]. In this paper, we give a new discovery of conversion of generalized singular values. New formulation is proposed.

\section{Preliminaries}In this section we give some definitions and notations.\\
$\mathbf{Definition 2.1}$
\emph{Let $A\in\mathbb{C}^{m\times n}$ and $B\in\mathbb{C}^{p\times n}$. A matrix pair $\{A,B\}$ is an $(m,p,n)$-Grassman matrix pair (GMP) if rank $(A^{H},B^{H})^{H}=n$. Let $A\in\mathbb{R}^{m\times n}$ and $B\in\mathbb{R}^{p\times n}$. A matrix pair $\{A,B\}$ is an $(m,p,n)$ real matrix pair (RMP) if rank $(A^{T},B^{T})^{T}=n$.}\\
$\mathbf{Definition 2.2}$
\emph{Let $\{A,B\}$ be an $(m,p,n)$-GMP (RMP). A nonnegative number-pair $(\alpha,\beta)$ is a generalized singular value of
the GMP (RMP) $\{A,B\}$ if
\begin{equation*}
(\alpha,\beta)=(\sqrt{\lambda},\sqrt{\mu}),\quad where\;(\lambda,\mu)\in\lambda(A^{H}A,B^{H}B)\; and\;\lambda,\mu\geq0.
\end{equation*}
The set of GSV of $\{A,B\}$ is denoted by $\sigma\{A,B\}.$ Evidently
\begin{equation*}
\sigma\{A,B\}=\{(\alpha,\beta)\neq(0,0)\;|\; \mathrm{det}(\beta^{2}A^{H}A-\alpha^{2}B^{H}B)=0,\;\alpha,\beta\geq0\}.
\end{equation*}
Let $\{A,B\}$ be an $(m,p,n)$-GMP and SVD of $\left(
\begin{array}{c}
A \\
B
\end{array}
\right)$ be
\begin{eqnarray}\label{sim}
\left(
\begin{array}{c}
A \\
B
\end{array}
\right)=L\Upsilon K^{H}=\left(
\begin{array}{cc}
L_{1}^{H}& L_{11}\\
L_{2}^{H}&L_{21}
\end{array}
\right)\left(
\begin{array}{c}
\Upsilon_{1} \\
O_{(m+p-n)\times n}
\end{array}
\right) K^{H},
\end{eqnarray}
where $L\in\mathbb{U}_{m+p},\Upsilon\in\mathbb{C}^{(m+p)\times n},\Upsilon_{1}\in\mathbb{C}^{n\times n}, K\in\mathbb{U}_{n},L_{1}^{H}\in\mathbb{C}^{m\times n}$. Let singular values of $L_{1}$ and $L_{2}$ be $\gamma_{i},\theta_{i}$ with $\gamma_{1}\geq\gamma_{2}\geq\cdots\geq\gamma_{\min\{m,n\}}\geq0$ and $0\leq\theta_{1}\leq\theta_{2}\leq\cdots\leq\theta_{\min\{p,n\}}$.}

\section{Conversion of GSVs of GMP by singular values}We give a simple conversion of generalized singular values of Grassmann matrix pairs by singular values.\\
$\mathbf{Theorem 3.1}$
Let $\{A, B\}$ be an $(m, p, n)$ Grassman matrix pair (real matrix pair) and $\sigma\{A,B\}=\{(\alpha_{i},\beta_{i})\}_{i=1}^{n}$. Then
for $n\leq m$,
$\alpha_{i}=\gamma_{i},\;i=1,\ldots,n$ and for $m\leq n$,
$\alpha_{i}=\gamma_{i},\;i=1,\ldots,m,\;\alpha_{m+1}=\cdots=\alpha_{n}=0.$
For $n\leq p$,
$\beta_{i}=\theta_{i},\;i=1,\ldots,n$ and for $p\leq n$,
$\beta_{n-p+i}=\theta_{i},\;i=1,\ldots,p,\;\beta_{1}=\cdots=\beta_{n-p}=0.$\\
Proof. It follows from Definition 2.1 that
\begin{equation*}
A(A^H A+B^H B)^{-\frac{1}{2}}=U\Sigma_{A}R(R^{H}R)^{-\frac{1}{2}},
\end{equation*}
\begin{equation*}
B(A^H A+B^H B)^{-\frac{1}{2}}=V\Sigma_{B}R(R^{H}R)^{-\frac{1}{2}}.
\end{equation*}
Since $U,V,R(R^{H}R)^{-\frac{1}{2}}$ are unitary, then singular values of $A(A^H A+B^H B)^{-\frac{1}{2}}$
are the same as singular values of $\Sigma_{A}$
and singular values of $B(A^H A+B^H B)^{-\frac{1}{2}}$ are the same as singular values of $\Sigma_{B}$. It follows that
\begin{eqnarray}\label{l1}
A(A^H A+B^H B)^{-\frac{1}{2}}&=&(L_{1}^{H},L_{11})\Upsilon K^{H}K(\Upsilon^{H}\Upsilon)^{-1/2}K^{H}\nonumber\\
&=&(L_{1}^{H},L_{11})\Upsilon (\Upsilon^{H}\Upsilon)^{-1/2}K^{H}\nonumber\\
&=&(L_{1}^{H},L_{11})\left(
\begin{array}{c}
\Upsilon_{1} \\
O_{(m+p-n)\times n}
\end{array}
\right) (\Upsilon_{1}^{H}\Upsilon_{1})^{-1/2}K^{H}=L_{1}^{H}K^{H}
\end{eqnarray}
\begin{eqnarray}\label{l2}
B(A^H A+B^H B)^{-\frac{1}{2}}&=&(L_{2}^{H},L_{21})\Upsilon K^{H}K(\Upsilon^{H}\Upsilon)^{-1/2}K^{H}\nonumber\\
&=&(L_{2}^{H},L_{21})\Upsilon (\Upsilon^{H}\Upsilon)^{-1/2}K^{H}\nonumber\\
&=&(L_{2}^{H},L_{21})\left(
\begin{array}{c}
\Upsilon_{1} \\
O_{(m+p-n)\times n}
\end{array}
\right) (\Upsilon_{1}^{H}\Upsilon_{1})^{-1/2}K^{H}=L_{2}^{H}K^{H}
\end{eqnarray}
where $K^{H}\in\mathbb{U}_n$. Hence, singular values of $A(A^H A+B^H B)^{-\frac{1}{2}}$
are the same as singular values of $L_{1}$
and singular values of $B(A^H A+B^H B)^{-\frac{1}{2}}$ are the same as singular values of $L_{2}$.
Then by Definition 2.1 and Definition 2.2 we have
for $n\leq m$,
$\alpha_{i}=\gamma_{i},\;i=1,\ldots,n$ and for $m\leq n$,
$\alpha_{i}=\gamma_{i},\;i=1,\ldots,m,\;\alpha_{m+1}=\cdots=\alpha_{n}=0.$
For $n\leq p$,
$\beta_{i}=\theta_{i},\;i=1,\ldots,n$ and for $p\leq n$,
$\beta_{n-p+i}=\theta_{i},\;i=1,\ldots,p,\;\beta_{1}=\cdots=\beta_{n-p}=0.$\\
$\mathbf{Remark 3.1}$
\emph{By Theorem 3.1 we have computing GSV can be converted to calculating singular values of $m\times n$ matrix $L_{1}$ and
$p\times n$ matrix $L_{2}.$ Meanwhile, we note $\alpha_{i}$ and $\beta_{i}$ satisfy $\alpha_{i}^2+\beta_{i}^2=1$, then we need solve one
of $\alpha_{i}$ and $\beta_{i}$ and the other one can be computed naturally. Therefore, if $m\leq p$ we compute singular values of $m\times n$ matrix $L_{1}$ to deduce $\alpha_{i}$, and $\beta_{i}$ is deduced naturally. If
$p\leq m$ we compute singular values of $p\times n$ matrix $L_{2}$ to deduce $\beta_{i}$, and $\alpha_{i}$ is deduced naturally.}

\section{References}
[1] Ian N. Zwaan, Michiel E. Hochstenbach, Generalized Davidson and multidirectional-type methods for the generalized singular value decomposition, arXiv:1705.06120.

[2] Z. Drma$\breve{\mathrm{c}}$, A tangent algorithm for computing the generalized singular value decomposition, SIAM J. Numer. Anal., 35(1998), pp. 1804-1832.

[3] S Friedland,  A new approach to generalized singular value decomposition, SIAM J. Matrix Anal. Appl.,27(2005), pp. 434-444.

\end{document}